\documentclass[9pt,twocolumn,twoside]{arXiv}

\templatetype{arXiv_researcharticle} 

\title{A Hamilton-Jacobi-based  Proximal Operator}

\author[a,1]{Stanley Osher}
\author[b,1]{Howard Heaton}
\author[c,1,2]{\large Samy Wu Fung}

\affil[a]{Dept. of Mathematics, University of California, Los Angeles}
\affil[b]{Typal Academy}
\affil[c]{Dept. of Applied Mathematics and Statistics, Colorado School of Mines}

\usepackage{lipsum}
\newcommand{\ie}{\textit{i.e.}\ }
\newcommand{\eg}{\textit{e.g.}\ }
\newcommand{\prox}[1]{\mathrm{prox}_{#1}}
\newcommand{\bbR}{\mathbb{R}}
\newcommand{\bbN}{\mathbb{N}}
\newcommand{\bbP}{\mathbb{P}}
\newcommand{\bbE}{\mathbb{E}}
\newcommand{\sO}{\mathcal{O}}
\newcommand{\sN}{\mathcal{N}}
\newcommand{\sS}{\mathcal{S}}
\newcommand{\sB}{\mathcal{B}}

\newcommand{\itemsymbol}{{\small $\blacktriangleright$}}

\usepackage{amsthm} 
\newtheorem{theorem}{Theorem} 

\newtheorem{lemma}{Lemma}

\newtheorem{remark}{Remark}
\newtheorem{definition}{Definition} 

\DeclareMathOperator*{\argmin}{argmin}
\usepackage{subfig}

\usepackage{algorithm,setspace}
\usepackage{algcompatible}

\leadauthor{Lead author last name} 

\significancestatement{\rmfamily 
Many objectives do not admit explicit proximal formulas (\eg when objectives are either nonconvex or only accessible via an oracle). Yet, only using (possibly noisy) objective samples, we give a formula for accurately approximating such proximals.}

\authorcontributions{\ }
\authordeclaration{Code is available at \href{https://hj-prox.research.typal.academy}{hj-prox.research.typal.academy}}
\equalauthors{\textsuperscript{1}SO (Stanley Osher),  HH (Howard Heaton), and SWF (Samy Wu Fung) contributed equally.}
\correspondingauthor{\textsuperscript{2}Correspondence should be addressed to SWF via email:  swufung@mines.edu.}


\begin{abstract} 
First-order optimization algorithms are widely used today. Two standard building blocks in these algorithms are proximal operators (proximals) and gradients. Although gradients can be computed for a wide array of functions, explicit proximal formulas are only known for limited classes of functions.  We provide an algorithm, HJ-Prox, for accurately approximating such proximals. This is derived from a collection of relations between proximals, Moreau envelopes, Hamilton-Jacobi (HJ) equations, heat equations, and importance sampling. In particular, HJ-Prox smoothly approximates the Moreau envelope and its gradient. The smoothness can be adjusted to act as a denoiser. Our approach applies even when functions are only accessible by (possibly noisy) blackbox samples. We show HJ-Prox is effective numerically via several examples.
\end{abstract}

\dates{\today} 
\begin{document}

\maketitle
\thispagestyle{firststyle}
\ifthenelse{\boolean{shortarticle}}{\ifthenelse{\boolean{singlecolumn}}{\abscontentformatted}{\abscontent}}{}

\dropcap{T}he rise of computational power and availability of big data brought great interest to first-order optimization methods.
Second-order methods (\eg Newton's method) are effective with moderately sized problems, but generally do not scale well due to     memory requirements increasing quadratically with problem size and computation costs increasing cubically. First-order methods are often comprised of gradient and proximal operations, which are typically cheap to evaluate relative to problem size. Although gradients can be computed for many functions (or numerically approximated), the computation of proximals involves solving a small optimization problem. In   special cases (\eg with $\ell_1$ norms), these subproblems admit closed-form solutions  that can be quickly evaluated (\eg see~\cite{beck2017first}). These formulas yield great utility in many applications. However, we are presently interested in the   class of problems with (potentially nondifferentiable) objectives for which \textit{proximal formulas are unavailable}. 

We propose a new approach to compute proximal operators and corresponding Moreau envelopes for   functions $f$. 
We leverage the fact that the Moreau envelope of $f$ is the solution to a Hamilton-Jacobi (HJ) equation~\cite{evans2010partial}. 
The core idea is to add artificial viscosity to HJ equations and obtain explicit formulas for the proximal and Moreau envelopes using Cole-Hopf  transformation~\cite[Sec. 4.5.2]{evans2010partial}.
This approach enables proximals and Moreau envelopes of arbitrary $f$ to be approximated. 
Our proposed proximal approximations (called HJ-Prox) are computed \emph{using only function evaluations} and can, thus, be used in a zeroth-order fashion when integrated within an optimization algorithm. Finally, an importance sampling procedure is employed to mitigate the curse of dimensionality when estimating the HJ-Prox in dimensions higher than three. Numerical experiments show HJ-Prox is effective when employed within optimization algorithms when the proximal is unavailable and for blackbox oracles. 
Our work can generally be applied to first-order proximal-based algorithms such as Alternating Direction Method of Multipliers (ADMM) and its variants~\cite{powell1969method, boyd2011distributed, hestenes1969multiplier, fung2020admm}, and operator splitting algorithms~\cite{eckstein1992douglas,lions1979splitting,ryu2022large, davis2017three, goldstein2009split}.

\section*{Proximal Operators and Moreau Envelopes}

\begin{figure}
    \centering 
    \includegraphics[width=0.4\textwidth]{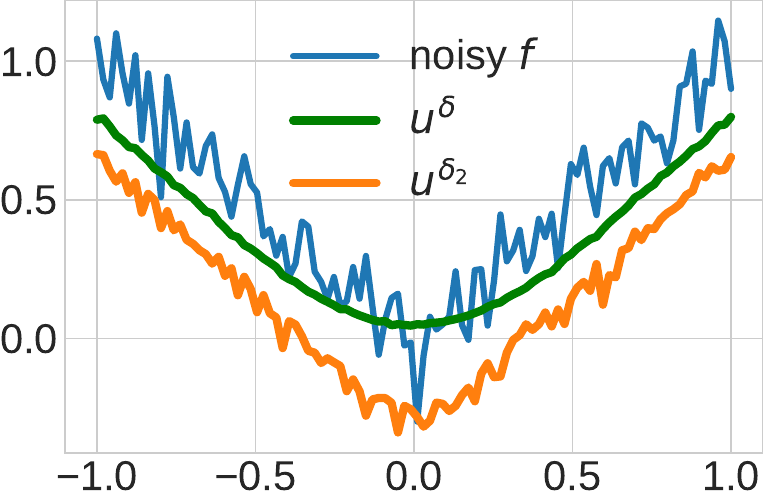}  
    \caption{Moreau envelope approximation $u^\delta$ using only noisy function samples with $\delta =0.1$ and $\delta_2 = 0.01$.}
    \label{fig: eye-candy} 
\end{figure}

Consider a function $f\colon\bbR^n\rightarrow\bbR$ and time $t > 0$.
The proximal $\prox{tf}$   and the Moreau envelope $u $ of $f$ \cite{moreau1962decomposition,bauschke2017convex} are defined by
\begin{equation}
    \prox{tf}(x) \triangleq \argmin_{z\in\bbR^n} f(z) + \dfrac{1}{2t}\|z-x\|^2
    \label{eq: original_prox_formula}
\end{equation}
and
\begin{equation} 
    u(x,t) \triangleq \min_{z\in \bbR^n} f(z) + \dfrac{1}{2t}\|z-x\|^2.
    \label{eq: moreau-envelope-definition}
\end{equation}
The proximal is the set of minimizers defining the envelope.
As shown in Figure \ref{fig: eye-candy}, the envelope $u$ widens valleys of $f$ while sharing global minimizers.  A well-known result (\eg see \cite{rockafellar1970convex,beck2017first}) states, if the envelope $u$ is differentiable at $x$, then 
\begin{align}
    \nabla u(x, t)= \dfrac{x-  \prox{t f}(x)}{t} .
\end{align}
Rearranging reveals
\begin{equation}
    \prox{tf}(x) = x - t\nabla u(x,t).
    \label{eq: prox-identity-grad-u}
\end{equation}
A key idea we use is to estimate the proximal by replacing $u$ with a smooth approximation $u^\delta \in C^\infty(\bbR)$, derived from a Hamilton-Jacobi (HJ) equation.

\section*{Hamilton-Jacobi Connection}
The envelope $u$ is a special case of the Hopf-Lax formula \cite{evans2010partial}. Fix any $T > 0$. For all $t\in [0, T]$, the envelope $u$ is a viscocity solution (\eg see \cite[Theorem 3.2]{evans2014envelopes})  to the HJ equation 
\begin{equation}
    \left \lbrace\begin{array}{rll}
    u_t  + \dfrac{1}{2}\|\nabla u \|^2 \hspace*{-7pt}& = 0  & \mbox{in\  $\bbR^n\times (0,T]$}\\
    u\hspace*{-7pt} & = f & \mbox{on\  $\bbR^n\times\{t=0\}$}.
    \end{array}\right.
    \label{eq: HJ}
\end{equation} 
Fixing $\delta > 0$, the associated viscous HJ equation is
\begin{equation}
     \left\lbrace \begin{array}{rll}
    u_t^\delta  + \frac{1}{2}\|\nabla u^\delta  \|^2 \hspace*{-7pt}& = \frac{\delta}{2} \Delta u^\delta   & \mbox{in\  $\bbR^n\times (0,T]$}\\[5pt]
    u^\delta \hspace*{-7pt} & = f  & \mbox{on\  $\bbR^n\times\{t=0\}$}.
    \end{array}  \right. 
    \label{eq: HJ-viscous}
\end{equation}
If $f$ is bounded and Lipschitz,   Crandall and Lions \cite{crandall1984two} show  $u^\delta$ approximates $u$, \ie $u^\delta \rightarrow u$ uniformly as $\delta \rightarrow 0^+$.  

\section*{Cole-Hopf Transformation}
Using the  transformation  $v^\delta \triangleq \exp(-u^\delta/\delta)$, originally attributed to Cole and Hopf~\cite{evans2010partial,chaudhari2018deep}, the function $v^\delta$ solves the heat equation, \ie 
\begin{equation}
    \left\lbrace \begin{array}{rll}
    v^\delta _t  -   {\frac{\delta}{2}} \Delta v^\delta \hspace*{-7pt}&= 0  & \mbox{in\  $\bbR^n\times (0,T]$}\\
    v^\delta \hspace*{-7pt} & = \exp(-f/\delta)  & \mbox{on\  $\bbR^n\times\{t=0\}$}.
    \end{array}\right.
    \label{eq: heat-PDE}    
\end{equation}
This transformation   is of   interest since $v^\delta$ can be expressed via the convolution formula (\eg see \cite{evans2010partial} for a derivation)
\begin{subequations}
    \begin{align}
    v^\delta(x,t) &=\Big( \Phi_{{\delta} t} * \exp(-f/\delta)\Big)(x)
    \\
    &= (2 \pi {\delta} t)^{-n/2} \int_{\bbR^n} \Phi_{\delta t}(x-y) \exp\left( -f(y) / \delta  \right) \mbox{d}y,
    \end{align}
\end{subequations}
where $\Phi_{{\delta t}}$ is a fundamental solution to [\ref{eq: heat-PDE}], \ie 
\begin{equation}
    \Phi_{{\delta} t}(x) \triangleq \left\lbrace \begin{array}{cl} ({2}\pi {\delta} t)^{-n/2}  \exp\left(-|x|^2/(2\delta t)  \right) & \mbox{in $\bbR^n\times (0,\infty)$}\\
    0 & \mbox{otherwise.}\end{array}\right.
\end{equation}
Using algebraic manipulations, we   recover the viscous   solution
\begin{equation}
    u^\delta(x,t) = - \delta \ln\Big(\Phi_{{\delta} t} * \exp(-f/\delta)\Big)(x)
    \ \ \mbox{in\  $\bbR^n\times (0,T]$}.
    \label{eq: u_delta_def}
\end{equation}
Differentiating reveals 
\begin{equation}
    \nabla u^\delta(x,t)
    = -\delta \cdot \nabla\left[   \ln\left(v^\delta(x,t)\right)\right] 
     = -\delta \cdot \dfrac{\nabla v^\delta(x,t)}{v^\delta(x,t)}.
     \label{eq: u-viscous-gradient}
\end{equation}

\section*{Importance Sampling}
\begin{figure*}[t]
    \centering
    \begin{tabular}{cccc} 
        \multicolumn{4}{c}{\large $f(x) = \|x\|_1$, $ \; \; \text{prox}_{tf}(x) = \text{shrink}_t(x)$}
        \\
        \textbf{a1)} $f$, $u$, and $u^\delta$ & \textbf{b1)} Proximal Comparison & \textbf{c1)} Proximal Err vs \# Samples & \textbf{d1)} $u^\delta$ with noisy Samples
        \\
        \includegraphics[width=0.23\textwidth]{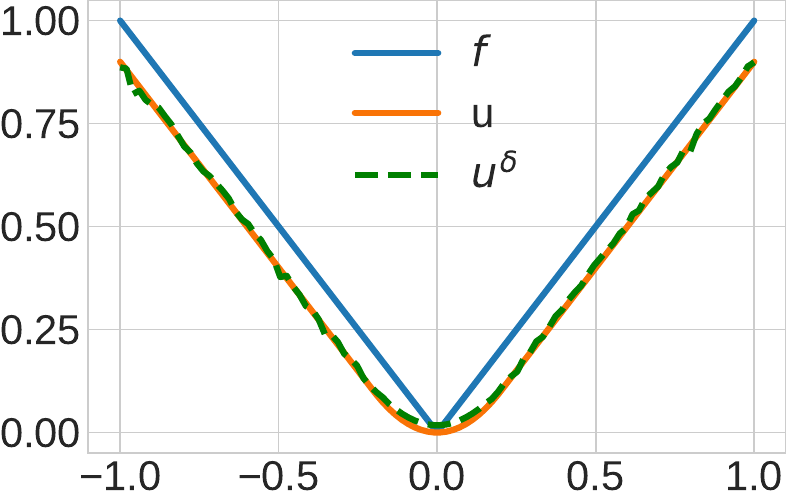}
        & 
        \includegraphics[width=0.23\textwidth]{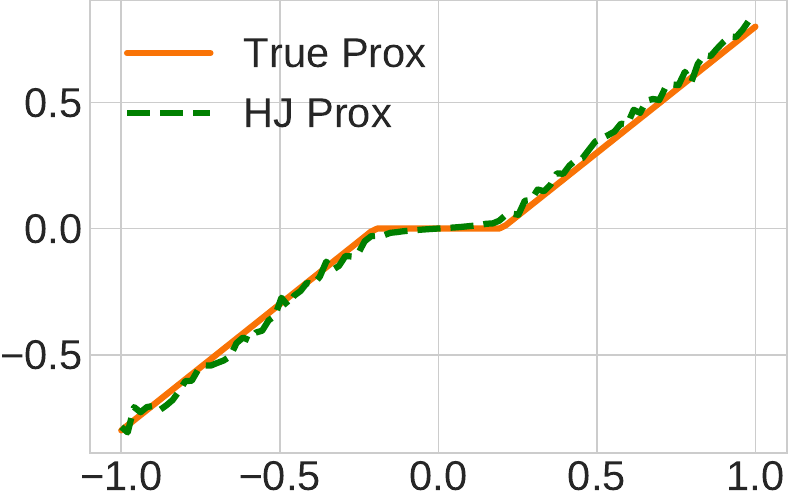}
        & 
        \includegraphics[width=0.23\textwidth]{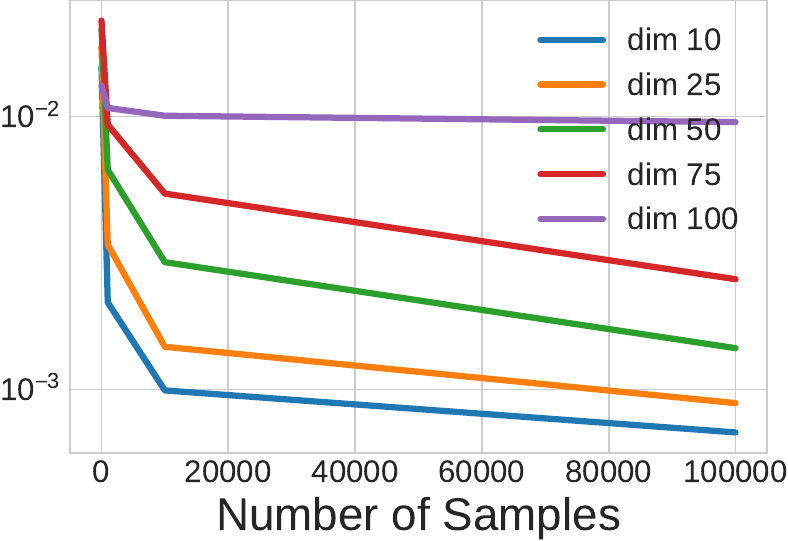}
        &
        \includegraphics[width=0.23\textwidth]{figures/l1_norm_envelope_noisy.pdf} 
        \\
        \\
        \multicolumn{4}{c}{\large $f(x) = \|x\|^2 + b^\top x$, $ \; \; \text{prox}_{tf}(x) = \frac{x-tb}{1+t}$, \;\; ($b = \mathbf{1}$)}
        \\
        \textbf{a2)} $f$, $u$, and $u^\delta$ & \textbf{b2)} Proximal Comparison & \textbf{c2)} $u^\delta$ from Noisy Samples & \textbf{d2)} Proximal from Noisy Samples
        \\
        \includegraphics[width=0.23\textwidth]{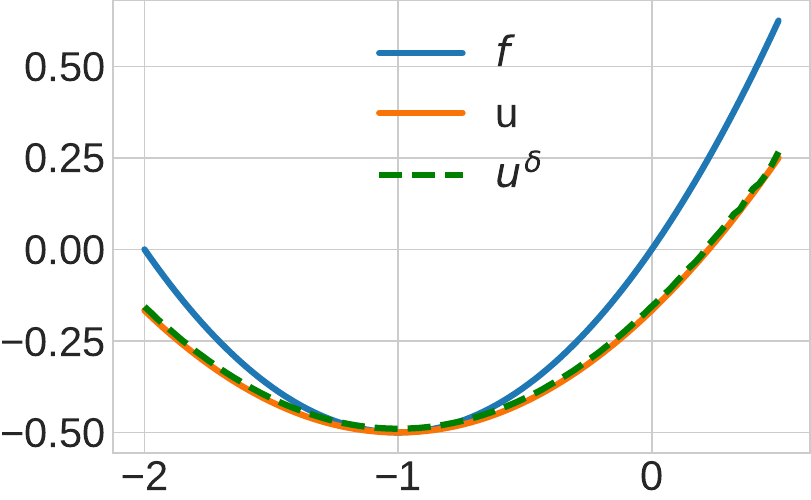}
        & 
        \includegraphics[width=0.23\textwidth]{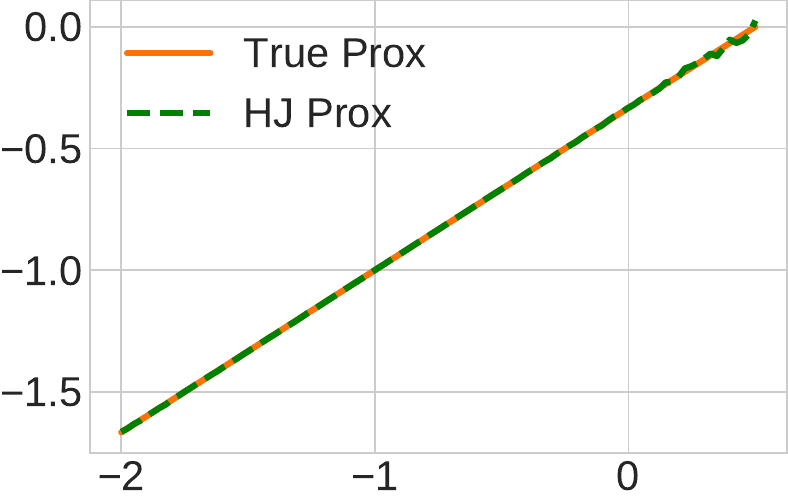}
        &
        \includegraphics[width=0.23\textwidth]{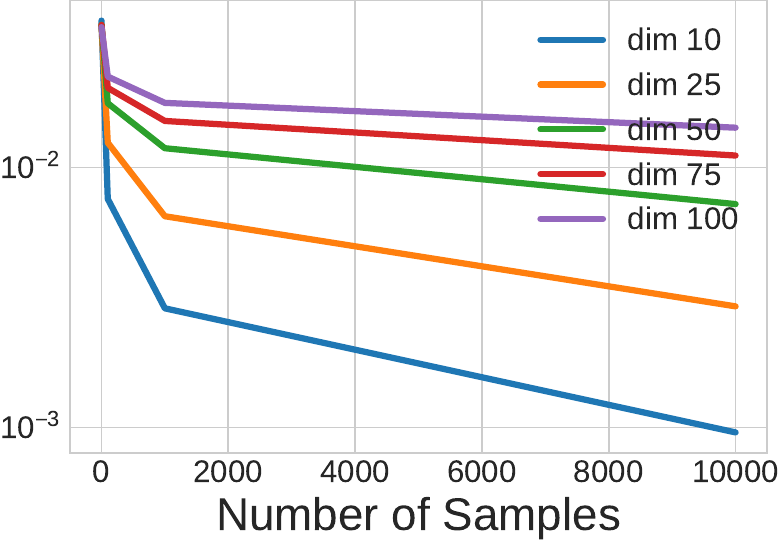}
        &
        \includegraphics[width=0.23\textwidth]{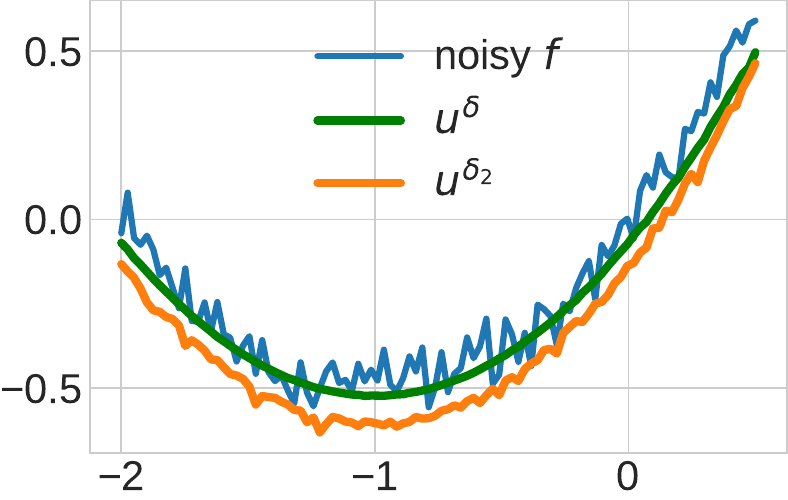}
        \\
        \\
        \multicolumn{4}{c}{\large $f(x) = -\sum_{i=1}^N \log(x_i)$, $ \; \; \text{prox}_{tf}(x)_i = \frac{x_i - \sqrt{x_i^2 + 4t}}{2}$}
        \\
        \textbf{a3)} $f$, $u$, and $u^\delta$ & \textbf{b3)} Proximal Comparison & \textbf{c3)} $u^\delta$ from Noisy Samples & \textbf{d3)} Proximal from Noisy Samples
        \\
        \includegraphics[width=0.23\textwidth]{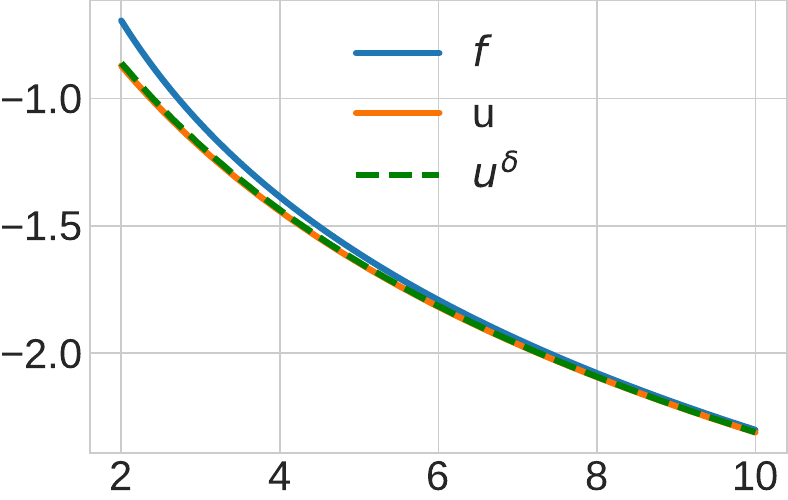}
        & 
        \includegraphics[width=0.23\textwidth]{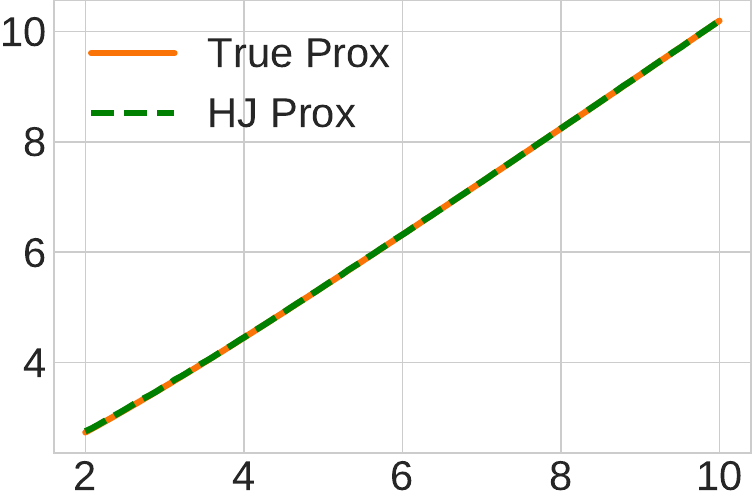}
        &
        \includegraphics[width=0.23\textwidth]{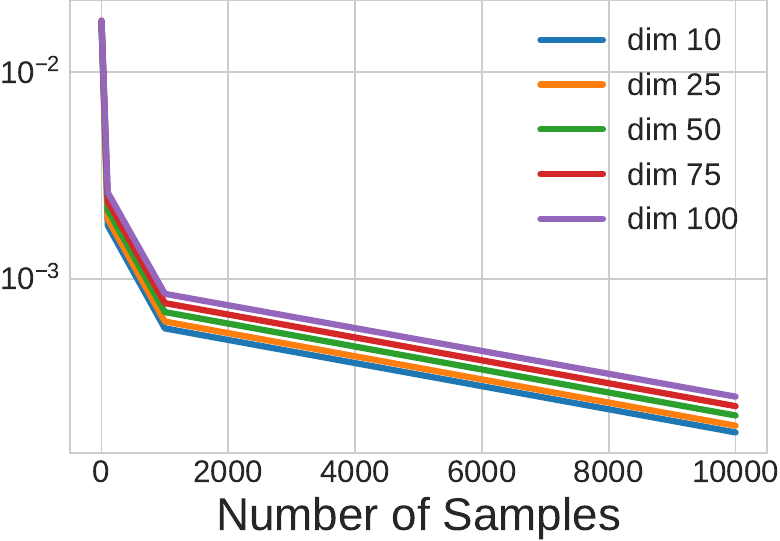}
        &
        \includegraphics[width=0.23\textwidth]{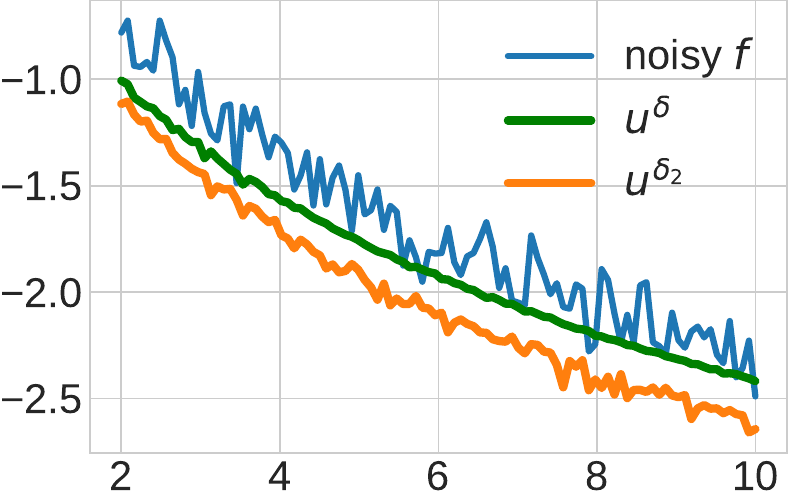} 
    \end{tabular}
    \caption{\textbf{(a1, a2, a3):} Plots for function $f$, exact Moreau envelope $u$, and HJ-based Moreau envelope $u^\delta$. \textbf{(b1, b2, b3):} Plots for true proximal and approximate HJ-based proximal operators. 
    \textbf{(c1, c2, c3):} Proximal approximations across different dimensions and samples.
    \textbf{(d1, d2, d3):} HJ-based Moreau envelopes $u^\delta$ obtained from \emph{noisy function samples}. Here, we use $\delta = 10^{-1}$ and $\delta_2 = 10^{-2}$. As expected, higher $\delta$ values have a stronger smoothing property. 
    The HJ-proximals are good approximations of the true proximal operators (seen through the Moreau envelopes) and can even be applied when only (potentially noisy) samples are available. 
    For the noisy case, we obtain a $C^\infty$ approximation of the underlying function $f$. For these experiments, we use $t = 0.1, 0.5, 2.0$ for rows 1, 2, and 3, respectively.}
    \label{fig: known_proximal_comparisons}
\end{figure*}

At first glance, the integral formula for $v^\delta$ in [\ref{eq: u-viscous-gradient}] may appear to require use of a grid for numerical estimation (and similarly for $\nabla v^\delta$). 
However, we  may avoid such grids by noting $v^\delta$ can be written as an expectation, \ie 
\begin{subequations}
\begin{align}
    v^\delta(x,t) 
    &= \Big( \Phi_{{\delta} t} * \exp(-f/\delta)\Big)(x)
    \\
    &= \bbE_{y\sim \sN(x,\delta t)}\left[ \exp\left(- f(y)/\delta \right) \right],
    \label{eq: v_delta}
\end{align}
\end{subequations}
where $y \sim \mathcal{N}(x, \delta t)$ denotes $y\in\bbR^n$ is sampled from a normal distribution with mean $x$ and standard deviation $\sqrt{\delta t}$.
In practice, finitely many samples $y^i \sim \mathcal{N}(x,\delta t)$  are used to estimate [\ref{eq: v_delta}]. This can greatly reduce sampling complexity~\cite{kloek1978bayesian, tokdar2010importance}.
Differentiating $v^\delta$ with respect to $x$ reveals
\begin{align}
    \nabla v^\delta(x,t)  
    = {-\frac{1}{\delta t}} \cdot \bbE_{y\sim  \sN(x,\delta t)}\left[(x-y) \exp\left(-{f}(y)/\delta \right) \right]
    \label{eq: nabla_v_delta}.
\end{align} 
\begin{algorithm}[t]
\newcommand{\alggap}{\hspace*{10pt}}
\caption{HJ-Prox -- Approximation of Proximal Operator}
\label{alg: HJ-prox}
\setstretch{1.35}
\begin{algorithmic}[1]  
    \STATE{HJ-Prox$(x,\ t;\ f,\  \delta,\  N, \ \alpha,\ \varepsilon):$} 
    \STATE{\alggap{\bf for} $i \in [N]$:}
    \STATE{\alggap \alggap Sample $y^i \sim \sN(x,\delta t/\alpha)$}
    \STATE{\alggap \alggap $z_i \leftarrow f(y^i)$}
    \STATE{\alggap \alggap \textbf{if} $z_i < 0$:}
    \STATE{\alggap \alggap \alggap {\bf return} HJ-Prox$(x,\ t;\ f + z_i + \varepsilon,\  \delta,\  N, \ \alpha,\ \varepsilon)$}
    \STATE{\alggap \alggap \textbf{if} $\exp\left(-\alpha z_i / \delta\right) \leq \varepsilon$:} 
    \STATE{\alggap \alggap \alggap {\bf return} HJ-Prox$(x,\ t;\ f,\  \delta,\  N, \ \alpha/2,\ \varepsilon)$}
 
    \STATE{\alggap prox $ \leftarrow \mbox{softmax}(-\alpha z/\delta)^\top [y^1 \cdots y^N]$}
    \STATE{\alggap {\bf return} prox}
\end{algorithmic} 
    \label{alg: recursive_hj_prox}
\end{algorithm} 
Plugging~[\ref{eq: v_delta}] and~[\ref{eq: nabla_v_delta}] into~[\ref{eq: u-viscous-gradient}] enables  $\nabla u^\delta$ to be written as
\begin{equation}
    \nabla u^\delta(x,t) = {\dfrac{1}{t}\cdot }\left({x} -  \dfrac{\bbE_{y\sim  \sN({x,{\delta} t})}\left[{y\cdot } \exp\left(-{f}(y)/\delta \right) \right]}
    {\bbE_{y\sim  \sN({x,{\delta} t})}\left[ \exp\left(-  {f}(y)/\delta \right) \right]}\right).
    \label{eq: u-viscous-gradient-expectation}
\end{equation}
The above relation was used in \cite{heaton2022global}.
Here we take a further step, combining [\ref{eq: prox-identity-grad-u}] and [\ref{eq: u-viscous-gradient-expectation}] to get an HJ-based estimate:
\begin{subequations}
    \begin{align}
        \prox{tf}(x)
        & = x - t \nabla u(x,t) \\
        & \approx x - t\nabla u^\delta(x,t) \\
        & =   \dfrac{\bbE_{y\sim  \sN({x,{\delta} t})}\left[y\hspace{0.5pt} \cdot  \exp\left(- {f}(y)/\delta \right) \right]}
    {\bbE_{y\sim  \sN({x,{\delta} t})}\left[ \exp\left(-  {f}(y)/\delta \right) \right]}.
    \end{align}\label{eq: prox_formula_expectation}\end{subequations} 
    As shown below, importance sampling enables efficient approximation of proximals in high dimensions (\eg see  Figure~\ref{fig: known_proximal_comparisons}). Moreover, [\ref{eq: prox_formula_expectation}] estimates proximals \emph{only using function values}, making it apt for zeroth-order optimization.

\section*{Numerical Considerations}
A possible numerical challenge in our formulation is to address numerical instabilities arising from the exponential term underflowing with limited numerical precision, due to either $\delta$ being small or $f(y)$ being large. To this end, note the proximal formula  may equivalently  be re-scaled via
\begin{subequations}
    \begin{align}
        \prox{tf}(x) &= \prox{\frac{t}{\alpha} \alpha f}(x)
        \label{eq: prox-alpha-equal}
        \\
        &\approx \dfrac{\bbE_{y\sim  \sN({x,{\delta} t/\alpha })}\left[{y\cdot } \exp\left(-\alpha {f}(y) / \delta \right) \right]}
        {\bbE_{y\sim  \sN({x,{\delta} t/\alpha})}\left[ \exp\left(-\alpha {f}(y)/\delta \right) \right]},
    \end{align}
\end{subequations}
where  $t$ is replaced by $t/\alpha$ and $f$ by $\alpha f$ in [\ref{eq: prox_formula_expectation}].
 
In this case, if $f/\delta$ becomes too large with respect to numerical precision limitations, it may be scaled down with a corresponding $\alpha$. To make the implementation stable, we check whether we obtain an underflow with $\exp(  \alpha f(y)/\delta)$ and rescale $\alpha$ using a linesearch-like approach. In particular, we recursively halve $\alpha$ until $\exp(  \alpha f(y)/\delta) > \varepsilon$ for a tolerance $\varepsilon$ (see line 7 of Algorithm~\ref{alg: recursive_hj_prox}. Yet,   small $\alpha$ makes the variance large and more samples may be required to accurately estimate the expectations. 
Another mitigation is to adaptively rescale $f$ based on the number of recursive steps taken in HJ-Prox.

Large $\delta$ can be used to smooth approximations and     mitigate   the stochastic characteristics of HJ-Prox.
Another potential instability that may arise is when $f$ is negative in certain parts of the domain. In this case, $\exp(  \alpha f(y)/\delta)$ may overflow. To remedy this, we   check whether $f(y)$ is negative and recursively shift the function until it is nonnegative (see line 5 of Algorithm~\ref{alg: recursive_hj_prox}).  

\section*{Convergence Analysis}
The   arguments above give intuition for a proximal approximation. However, having now the formula [\ref{eq: prox_formula_expectation}], we may formalize its utility without reference to differential equations.
Below we define two standard classes of functions used in optimization. 

\begin{definition}[Weakly Convex]
    For $\rho > 0$, a  function $f\colon\bbR^n\rightarrow\bbR$ is $\rho$-weakly convex if $f(x) + \frac{\rho}{2}\|x\|^2$ is convex.
\end{definition}

\begin{definition}[$L$-Smooth]
    For $L>0$, a function $f\colon\bbR^n\rightarrow\bbR$ is $L$-smooth if its gradient $\nabla f$ exists and is $L$-Lipschitz.
\end{definition}

\noindent Our main result shows HJ-Prox converges to the proximal.

\def\mainresult{If   $f\colon\bbR^n\rightarrow\bbR$ is $\rho$-weakly convex, for some $\rho>0$, and either $L$-smooth   or   $L$-Lipschitz, then, for all $x\in\bbR^n$, and $t\in (0,1/\rho)$, the proximal $\prox{tf}(x)$ is unique and, if $u(x,t) \geq  0$, then 
\begin{equation}
    \lim_{\delta \rightarrow 0^+} \dfrac{\bbE_{y\sim  \sN(x,\delta t/ )}\left[{y\cdot } \exp\left(- {f}(y) / \delta \right) \right]}
    {\bbE_{y\sim  \sN(x,\delta t/) }\left[ \exp\left(- {f}(y)/\delta \right) \right]}  = \prox{tf}(x).
\end{equation}}

\begin{theorem}[Proximal Approximation] \label{thm: proximal-approximation}
    \mainresult
\end{theorem}

\begin{remark}[Smoothing Property]
    In practice, we must pick positive $\delta$. Thankfully, increasing $\delta$ comes with the benefit of smoothing estimates (due to the Laplacian in the viscous HJ equation), as shown in rightmost column of Figure \ref{fig: known_proximal_comparisons}.
\end{remark}

\section*{Related Works}
Our proposal   closely relates to zeroth-order optimization algorithms, which do not require gradients. In fact, HJ-Prox does not require differentiability of $f$. 
Related methods include Random Gradients~\cite{ermoliev1988numerical, kozak2019stochastic, kozak2021stochastic,kozak2021zeroth}, sparsity-based methods~\cite{cai2022zeroth,cai2021zeroth, slavin2022adapting}, derivative-free quasi-Newton methods~\cite{berahas2019derivative,larson2019derivative, more2009benchmarking}, finite-difference-based methods~\cite{shi2021numerical, shi2021adaptive}, numerical quadrature-based methods~\cite{kim2021curvature, almeida1990learning}, Bayesian methods~\cite{larson2019derivative}, and comparison methods~\cite{cai2022one}. 
As proximals closely relate to gradient of Moreau envelopes, our work relates to methods that minimize Moreau envelopes (or their approximations)~\cite{chaudhari2019entropy,chaudhari2018deep, heaton2022global, scaman2020simple, scaman2020simple, davis2018stochastic, heaton2022global, doi:10.1137/18M1178244, davis2022escaping}.

The theoretical results in our work is closely related to the study of asymptotics as $\delta \to 0$ of integrals containing expressions of the form $\exp({-f/\delta})$, \ie  Laplace's method~\cite{evans2010partial}. 
Moreover, the idea of adding artificial diffusion to Burgers' equation and then applying Cole-Hopf transformation to approximate the gradient of the solution to the HJ equation has been largely developed in~\cite{evans2010partial} in the context of obtaining solutions to conservation laws in 1D.
The connections between Hopf-Lax and Cole-Hopf have been observed in the context of machine learning in~\cite{chaudhari2018deep}, image denoising and Bayesian inference~\cite{darbon2021bayesian, darbon2021connecting, louchet2008total}, and in the context of global optimization in~\cite{heaton2022global}.

\begin{figure*}
    \centering
    \begin{tabular}{cc}
        \multicolumn{2}{c}{\large \textbf{Moreau Envelope for Nonconvex Functions}}
        \\
        \\
        \large $f(x) = -\|x\|_1$, $ \; \; u(x) = -|x| - t/2$
        &
        \large $f(x) = -\frac{ax^2}{2}$, $ \; \; u(x) = -\frac{ax^2}{2(1-at)}, \quad t < \frac{1}{a}$
        \\
        \includegraphics[width=0.26\textwidth]{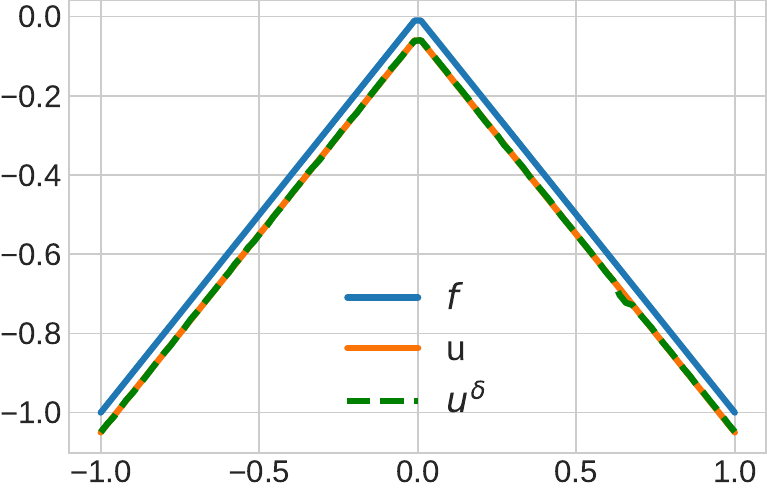}
        &
        \includegraphics[width=0.26\textwidth]{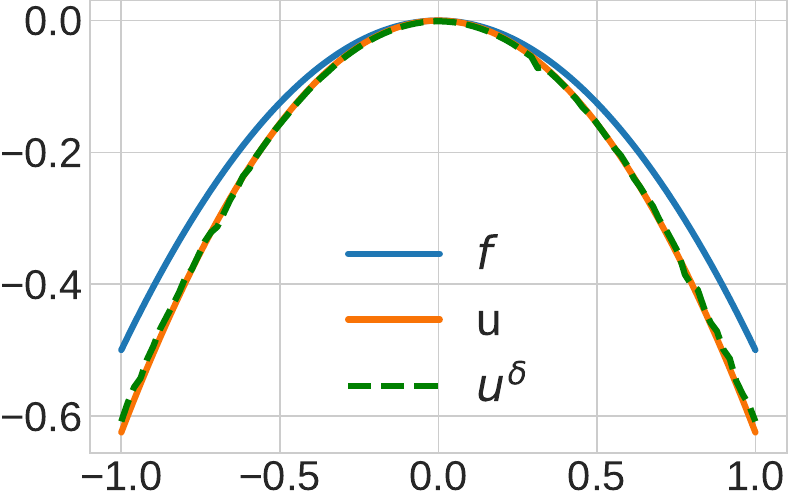}
    \end{tabular}
    \caption{HJ-based Moreau envelope for nonconvex functions with $t=0.1$ and $t=0.2$ in the left and right figures, respectively.}
    \label{fig: nonconvex_moreau_envelope_examples}
\end{figure*}

\begin{figure*}[t]
    \centering
    \begin{tabular}{cccc}
        \multicolumn{4}{c}{\large \textbf{Proximal Comparisons for Functions with Unknown Proximals}}
        \\
        \\
        \multicolumn{4}{c}{\large $f(x) = x^2 - \log(x)$}
        \\
        \textbf{a)} $f$, $u$, and $u^\delta$ & \textbf{b)} Proximal Comparison & \textbf{c)} $u^\delta$ from Noisy Samples & \textbf{d)} Proximal from Noisy Samples
        \\
        \includegraphics[width=0.23\textwidth]{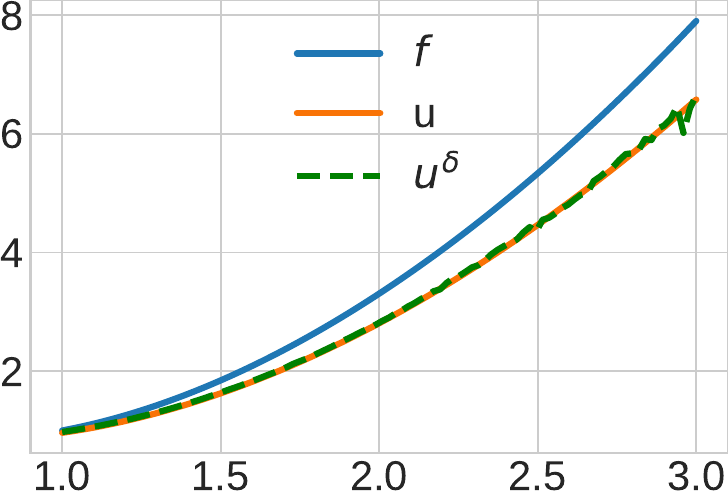}
        & 
        \includegraphics[width=0.23\textwidth]{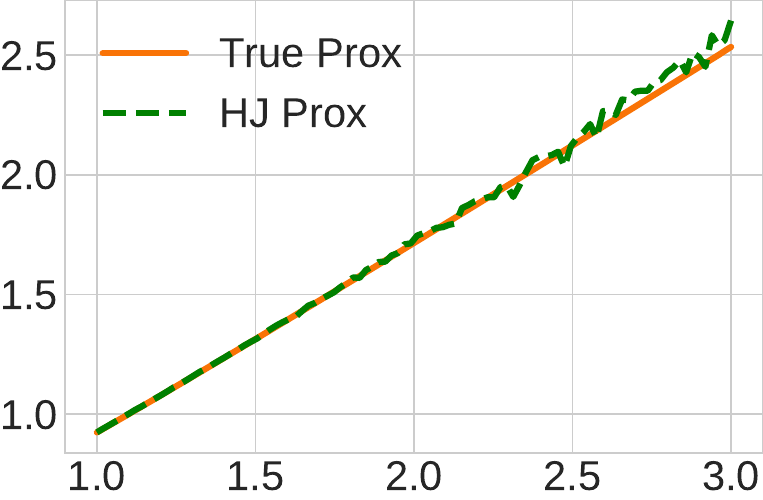}
        &
        \includegraphics[width=0.23\textwidth]{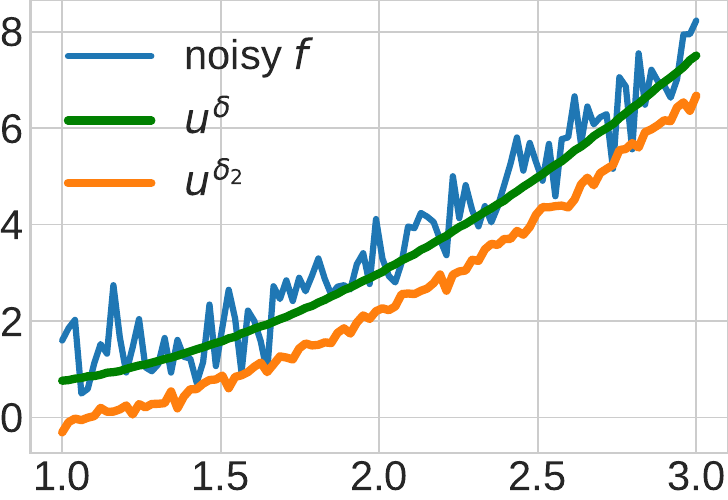}
        &
        \includegraphics[width=0.23\textwidth]{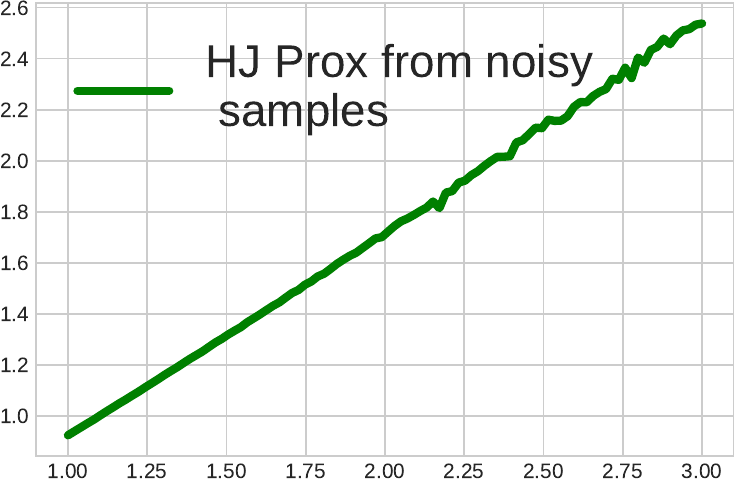}
    \end{tabular}
    \caption{\textbf{(a):} Plots for function $f$, exact Moreau envelope $u$, and HJ-based Moreau envelope $u^\delta$. \textbf{(b):} Plots for true proximal and approximate HJ-based proximal operators. \textbf{(c):} HJ-based Moreau envelopes $u^\delta$ obtained from \emph{noisy function samples}. \textbf{(d):} HJ-based proximal computed \emph{using noisy function samples}.
    Since there is no analytic proximal formula, we obtain the ``true'' proximal by solving the optimization~\eqref{eq: original_prox_formula} using gradient descent. The HJ-based proximal is a good approximation of the true proximal operators and can even be applied when only (potentially noisy) samples are available. As in the analytic case, we obtain a $C^\infty$ approximation of the underlying function $f$ in the noisy case. Here, $\delta = 0.1$ for the noiseless case and $\delta = 0.5$ and $\delta_2 = 0.1$ for the noisy case.}
    \label{fig: unknown_prox_function}
\end{figure*}

\begin{figure*}
    \centering
    \begin{tabular}{cc}
        \multicolumn{2}{c}{\large \textbf{HJ-ISTA Comparison}}
        \\
        \includegraphics[width=0.4\textwidth]{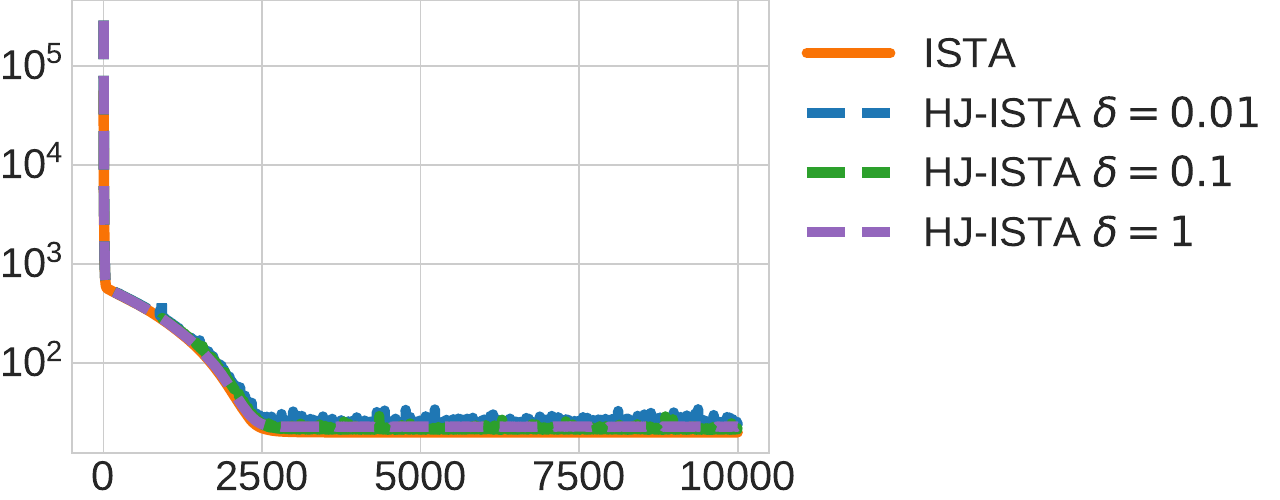}
        &
        \includegraphics[width=0.4\textwidth]{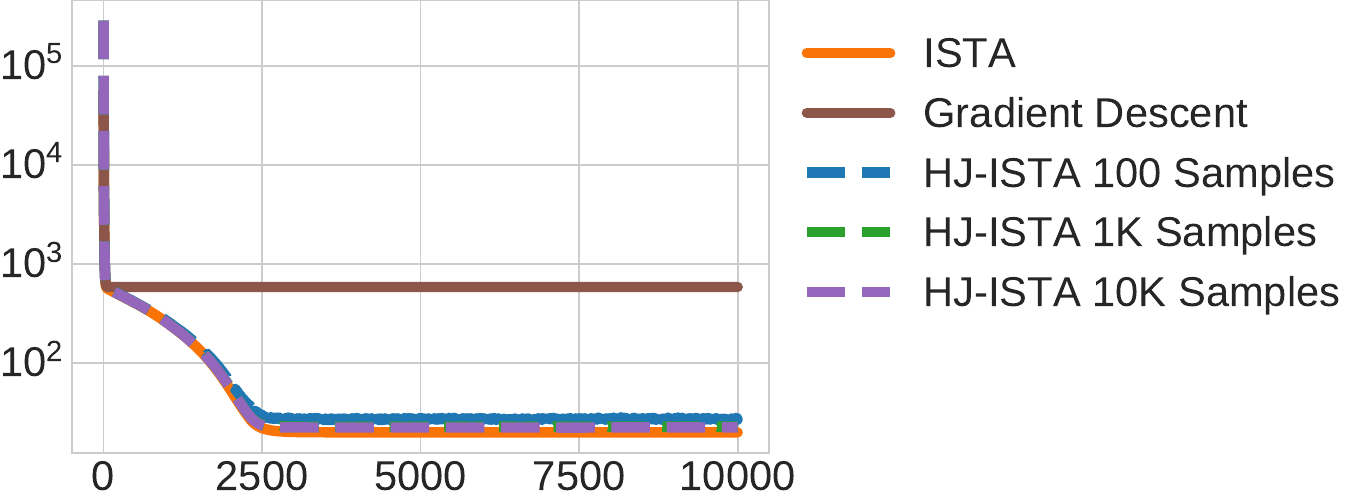}
        \\
        \textbf{(a)} Varying Smoothing $\delta$, Fixed \# Samples $N=1000$
        &
        \textbf{(b)} Varying \# Samples $N$, Fixed Smoothing $\delta=1.0$
    \end{tabular}
    \caption{Convergence plots showing function value for solution estimates $\{x^k\}$ when solving the LASSO problem [\ref{eq: lasso}] with ISTA, juxtaposing use of an analytic proximal formula, gradient descent (\ie  ignoring the proximal), and the approximate HJ-prox (Algorithm \ref{alg: HJ-prox}). 
    Plots with HJ-prox show averaged results from 30 trials with distinct random seeds. 
    To ensure the proximal is playing a role in the optimization process, we also show a function value history of gradient descent applied to the unregularized least squares problem in~[\ref{eq: lasso}] (\ie, with no $\ell_1$ norm term).
    }
    \label{fig: ista_comparison}
\end{figure*}

\section*{Numerical Experiments}

Examples herein show   HJ-Prox (Algorithm \ref{alg: HJ-prox}) can
\begin{itemize}[label=\itemsymbol]
    \item approximate proximals \textit{and} smooth noisy samples,
    \item converge comparably to existing algorithms, and
    \item solve a new class of zeroth-order optimization problems.
\end{itemize}
Each   item is addressed by a set of experiments. 
Regarding the last item, to our knowledge, HJ-Prox is the first tool to enable faithful solution estimation for constrained problems where the objective is only accessible via noisy blackbox samples. 

\subsection*{Proximal and Moreau Envelope Estimation}
Herein we compare HJ-Prox to known proximal operators. Figure~\ref{fig: known_proximal_comparisons} shows HJ-Prox for three functions (absolute value, quadratic, and log barrier) whose proximals are known. In the leftmost column (a), we show   the Moreau envelope $u(x,t)$ given by~[\ref{eq: moreau-envelope-definition}], and an estimate of Moreau envelope using the HJ-Prox $u^\delta(x,t)$. Given the close connection between proximals and Moreau envelopes, we believe this visual is a natural and intuitive way to gauge whether the proximal operator is accurate. 
Column (b) juxtaposes the true proximal and HJ-Prox. Column (c) shows the accuracy of HJ-Prox across different dimensions and numbers of samples. In the rightmost column (d), we estimate Moreau envelopes using HJ-Prox  \emph{using noisy function values}. The resulting envelopes are smooth since $u^\delta$ is a smooth (\ie $C^\infty$) approximation of $u$. Thus, HJ-Prox can be used to obtain smooth estimates from noisy observations.

Figure \ref{fig: nonconvex_moreau_envelope_examples} shows Moreau envelopes for nonconvex functions $f$. As in the other example, here HJ-based Moreau envelope estimates also accurately approximate Moreau envelopes.
Note these proximals may be well-defined only for small time $t$
(as the proximal operator objective in [\ref{eq: original_prox_formula}] is strongly convex for small $t$).  
 Lastly, we apply HJ-Prox with a function that has \emph{no analytic formula} for its proximal or Moreau envelope in Figure~\ref{fig: unknown_prox_function}. In this experiment, we obtain a ``true'' Moreau envelope and proximal operator by solving the minimization problem~[\ref{eq: original_prox_formula}] iteratively via gradient descent. Faithful recovery is shown in Figures~\ref{fig: unknown_prox_function}a and \ref{fig: unknown_prox_function}b, and smoothing in Figure \ref{fig: unknown_prox_function}c.

\begin{figure*}
    \centering
    \begin{tabular}{cc}
        \multicolumn{2}{c}{\large \textbf{Relative Errors for HJ-MM using noisy $f$}}
        \\
        \includegraphics[width=0.4\textwidth]{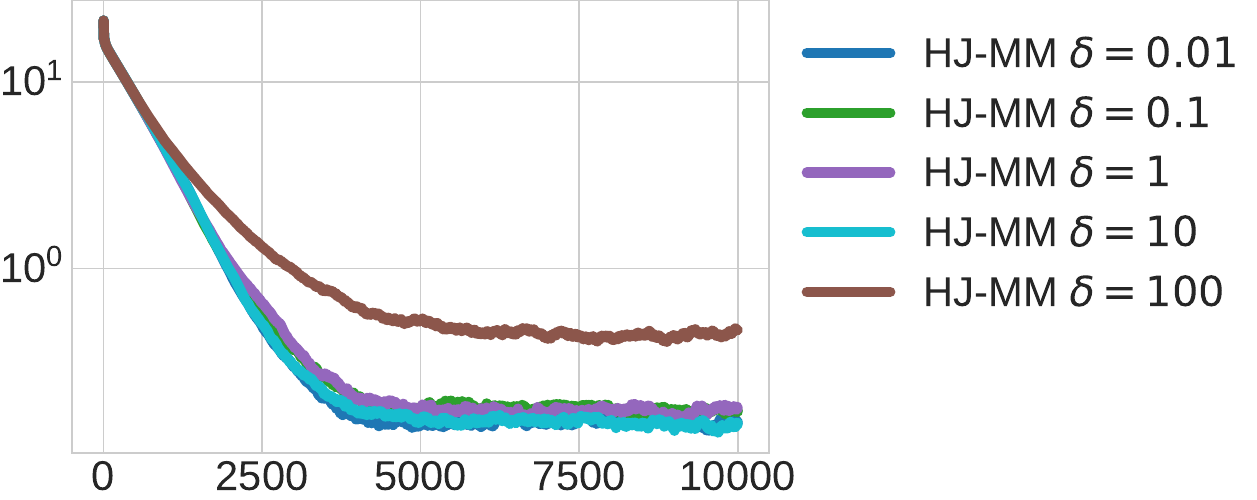}
        &
        \includegraphics[width=0.4\textwidth]{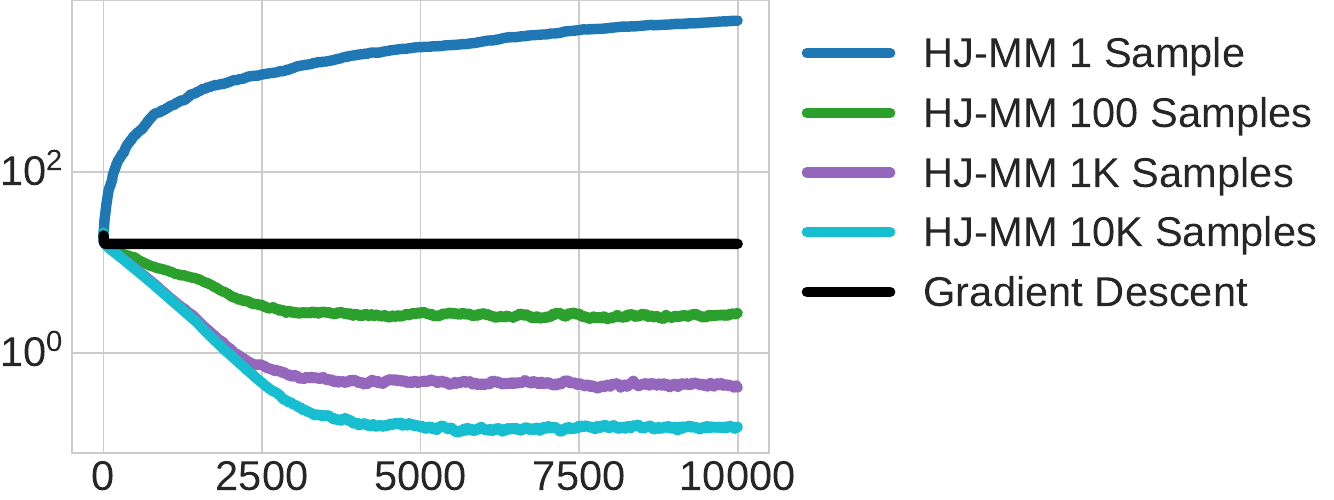}
        \\
        \textbf{(a)} Varying Smoothing $\delta$, Fixed \# Samples $N=$10K
        &
        \textbf{(b)} Varying \# Samples $N$, Fixed Smoothing $\delta=10$
    \end{tabular}
    \caption{Convergence plots showing relative errors for solution estimates $\{x^k\}$ when solving the minimization problem [\ref{eq: noisy_problem}] with linearized method of multipliers and HJ-prox (Algorithm 1). Each plot shows averaged results from 30 trials with distinct random seeds. 
    Due to the noise , we observe in (a) that a larger $\delta=10$ leads to a better approximation, but too large ($\delta=100$) leads to oversmoothing and   reduces accuracy. We find $\delta=10$ to be most optimal, and (b) shows that more samples lead to more accurate approximations (as expected). To ensure the proximal is playing a role in the optimization process, we also show the relative error when gradient decent is applied to the constraint residual in~[\ref{eq: noisy_problem}] (\ie, we only minimize constraint residual). Indeed, gradient descent performs poorly by comparison.}
    \label{fig: lmm-plots}
\end{figure*}

\subsection*{Optimization with Proximable Function}
This experiment juxtaposes HJ-prox and an analytic proximal formula in an optimization algorithm. Consider
 the Lasso problem~\cite{daubechies2004iterative} 
\begin{equation}
    \min_{x\in\bbR^{1000}} \frac{1}{2} \|Ax - b\|_2^2 + 0.1 \|x\|_1,
    \label{eq: lasso}
\end{equation}
where entries of $A \in \bbR^{500 \times 1000}$ and $b\in\bbR^{500}$ are i.i.d. Gaussian samples. The iterative soft thresholding algorithm (ISTA)~\cite{beck2009fast} defines a sequence of solution estimates $\{x^k\}$ for all $k\in\bbN$ via
\begin{equation}
    x^{k+1} = \mbox{shrink}\left( x^k - \beta A^\top(Ax^k - b);\ 0.01\beta \right), 
\end{equation}
where the shrink operator defined element-wise by
\begin{equation}
    \mbox{shrink}(x;\ t)
    \triangleq 
    \mbox{sign}(x) \max(0,\ |x| - t). 
    \label{eq: def-shrink}
\end{equation}
Figure~\ref{fig: ista_comparison} compares the convergence of ISTA using the shrink operator in [\ref{eq: def-shrink}] and HJ-Prox estimates of the shrink. To ensure convergence, we choose $\beta =1/\|A^\top A\|_2$. Our experiments show HJ-based ISTA can solve Lasso, up to an error tolerance.

\subsection*{Optimization with Noisy Objective Oracles}
Consider a constrained minimization problem where objective values $f$ can only be accessed via a noisy oracle\footnote{Here $\sO$ is a noisy function, \textit{not} to be confused with ``Big O'' often used to describe limit behaviors.} $\sO$. 
Our   task is to solve
\begin{equation} 
        \min_{x\in\bbR^{1000}}  \bbE[\sO(x)]
        \quad \mbox{s.t.}\quad Ax=b,
    \label{eq: noisy_problem}
\end{equation}
where $A$ and $b$ are as in the prior experiment and the expectation $\bbE$ is over  oracle noise.
To model ``difficult'' settings (\eg when a singular value decomposition of $A$ is unavailable), we do \textit{not}  use any projections onto the feasible set. As knowledge of the structure of $\sO$ is \textit{unknown} to the solver, we emphasize schemes for solving [\ref{eq: noisy_problem}] must   use zeroth-order optimization schemes \cite{larson2019derivative}. 
Here, each oracle call returns 
\begin{equation}
    \sO(x) = (1+\varepsilon)\cdot\|Wx\|_1,
    \quad \mbox{where}\ \varepsilon\sim \mathcal{N}(0, \sigma^2),
    \label{eq: noisy-oracle}
\end{equation}
with a \emph{new noise sample}  $\varepsilon \in \bbR$  used in each oracle evaluation, $\sigma=0.005$, and $W \in\bbR^{1000\times 1000}$  a fixed Gaussian matrix. In words, the noise has magnitude $0.5\%$ of $\|Wx\|_1$.
Although the oracle structure is shown by [\ref{eq: noisy-oracle}], our task is to solve [\ref{eq: noisy_problem}] \textit{without} such knowledge.  We   do this  with the linearized method of multipliers (\eg see Section 3.5 in \cite{ryu2022large}). Specifically, for each index $k\in\bbN$, the update formulas for the  solution estimates $\{x^k\}$ and  corresponding  dual variables $\{u^k\}$ are 
\begin{subequations}
    \begin{align}
        x^{k+1} & = \prox{t\sO}\left(x^k - tA^\top(u^k + \lambda (Ax^k - b))\right)\\
        u^{k+1} & = u^k + \lambda (Ax^{k+1} - b),
    \end{align}
\end{subequations} 
with step sizes $t = 1/\|A^\top A\|_2$ and $\lambda = 1/2$.
Without noise $\varepsilon$, convergence occurs if $t\lambda \|A^\top A\|_2 < 1$~\cite{ryu2022large}, justifying   our choices for $t$ and $\lambda$. 
The proximal $\prox{t\sO}$ is estimated by HJ-prox. 

We \textit{separately} solve the optimization problem using full knowledge of the objective $\|Wx\|_1$ without noise; doing this enables us to plot the relative error of the sequence $\{x^k\}$ in Figure \ref{fig: lmm-plots}. 
All the plots show $\{x^k\}$ converges to the optimal $x^\star$, up to an error threshold, regardless of the choice of $\delta$ and number of samples $N$.  
Notice Figure \ref{fig: lmm-plots}a shows ``small'' values of $\delta$ give comparable accuracy, but that oversmoothing with ``large'' $\delta=100$ degrades performance of the algorithm.
These plots also illustrate the HJ-prox formula is efficient with respect to calls to the oracle $\sO$. Indeed, note the plots in Figure \ref{fig: lmm-plots}b that decrease relative error use, at each iteration, respectively use 0.1, 1, and 10 oracle calls per dimension of the problem! We hypothesize the smoothing effect of the viscous $u^\delta$ and averaging effect of importance sampling contribute to the observed convergence. In   this experiment,  HJ-prox   converges to within an error tolerance, is efficient with respect to oracle calls, and smooths Gaussian noise.

\section*{Conclusion}
We propose a novel algorithm, HJ-prox, for efficiently approximating proximal operators. 
This is derived from approximating   Moreau envelopes via   viscocity solutions to Hamilton-Jacobi (HJ) equations, as given via the Hopf-Lax formula. 
Upon rewriting this approximation in terms of expectations, we use importance sampling to avoid discretizing the integrals, thereby mitigating the curse of dimensionality. Our numerical examples show HJ-Prox is effective for a collection of functions, both with and without known proximal formulas.  Moreover, HJ-prox can be effectively used in constrained optimization problems  \emph{even when only noisy objective values are available}.


\section*{Acknowledgements}
SO thanks the funding from AFOSR MURI FA9550-18-1-0502, ONR:N00014-20-1-2093
and N00014-20-1-2787, and NSF DMS 2208272 and 1952339.

\section*{References }
\bibliography{hj-prox-refs}

\newpage
\appendix
\section*{Proofs}

\normalsize 

For concise expression below, for $t > 0$ and $\delta > 0$ we define
\begin{equation}
    \phi_t(z) \triangleq f(z) + \dfrac{1}{2t}\|z-x\|^2,
\end{equation}
$\phi_t^\star \triangleq \inf\{\phi_t(y) : y\in\bbR^n\}$,
and
\begin{equation}
    \sigma_\delta (z)
        \triangleq
        \dfrac{\exp\left(-\phi(z)/\delta\right)}{\|\exp(-\phi/\delta)\|_{L^1(\bbR^n)}}.
        \label{eq: sigma-definition}
\end{equation}

\begin{lemma}\label{lemma: sigma-polynomial-vanish}
    If the conditions of Theorem \ref{thm: proximal-approximation} hold, then
    \begin{equation}
        \int_{\bbR^n} \sigma_\delta(y)\ \mathrm{d}y
        = 1,
        \quad 
        \sigma_\delta(y) \geq 0,\quad \mathrm{for\ all\ } y\in\bbR^n,
        \label{eq: proof-sigma-pdf}
    \end{equation}       
    and if   $r \in (0,1)$, then for all polynomials $p$ of positive degree
    \begin{equation}
        \lim_{\delta\rightarrow0^+}\int_{\bbR^n- \sB(\xi^\star,r)}
        \sigma_\delta(y) p(\|y-\xi^\star\|)\ \mbox{d}y = 0.
        \label{eq: sigma-polynomial-limit-vanish}
    \end{equation}
\end{lemma}

\begin{proof} 
    By algebraic limit laws, it suffices to verify [\ref{eq: sigma-polynomial-limit-vanish}] for any $p(x)=x^k$ with $k\geq 1$, and we proceed as follows. First we show $\sigma_\delta$ satisfies properties to be a probability density (Step 1). We show various $L^p$ norm limits hold for the numerator (Step 2) and denominator (Step 3) of integrating [\ref{eq: sigma-polynomial-limit-vanish}]. 
    Combining these limits gives [\ref{eq: sigma-polynomial-limit-vanish}] (Step 4).

    \paragraph{Step 1}
    The numerator and denominator in the definition [\ref{eq: sigma-definition}] for $\sigma_\delta$ are nonnegative, making $\sigma_\delta \geq 0$ everywhere.
    By the choice of $t$, $\phi_t$ is $\theta\triangleq 1/t - \rho$ strongly convex, and so it admits a unique minimizer $\xi^\star = \prox{tf}(x)$ and satisfies
    \begin{equation}
        \phi_t(y) \geq \phi_t^\star + \left<0, y-\xi^\star\right> + \dfrac{\theta}{2}\|y-\xi^\star\|^2,
        \quad \mbox{for all } y\in\bbR^n.
        \label{eq: proof-strong-convex}
    \end{equation}    Consequently, 
    \begin{equation}
        0 < e^{-\frac{\phi_t(y)}{\delta}}\leq  e^{-\frac{\phi_t^\star + \frac{\theta}{2}\|y-\xi^\star\|^2}{\delta}},
        \quad \mbox{for all } y\in\bbR^n.
        \label{eq: exponential-quadratic-bound}
    \end{equation}
    Since the upper bound above is an exponential that decays quadratically, the middle term in [\ref{eq: exponential-quadratic-bound}] is integrable over $\bbR^n$. 
    As $\phi_t^\star = u(\xi^\star, t)\geq 0$ by hypothesis,     the denominator in the definition of $\sigma_\delta$ is positive.
    Then [\ref{eq: proof-sigma-pdf}] readily follows.

    \paragraph{Step 2}  
    A classic  result in   analysis (\eg see  \cite[Exercise 3.4]{rudin1966real}) states $L^p$ norms converge to the $L^\infty$ norm as $p\rightarrow\infty$, and so
    \begin{equation}
        \lim_{\delta\rightarrow 0^+}
        \left\|e^{-\phi_t}\right\|_{L^{\frac{1}{\delta}}(\bbR^n)}
        = \|e^{-\phi_t}\|_{L^\infty(\bbR^n)}
        = e^{-\phi_t^\star},
        \label{eq: sigma-denominator-delta-limit}
    \end{equation}
    where  the $L^{\frac{1}{\delta}}$ norm is always finite by Step 1 and the final equality holds since the exponential is maximized by $\phi_t^\star$.

    \paragraph{Step 3} Integrating the numerator of [\ref{eq: sigma-polynomial-limit-vanish}] for $p(x)=x^k$ gives
    \begin{subequations} 
    \begin{align}
        &\int_{\bbR^n-\sB(\xi^\star,r)} e^{-\frac{\phi_t(y)}{\delta}} \|y-\xi^\star\|^k\ \mbox{d}y \\
        \leq &   \int_r^\infty e^{-\frac{\phi_t^\star + \frac{\theta \tau^2}{2}}{\delta}} \tau^{k}\cdot n|\sB(\xi^\star,1)|\tau^{n-1}\ \mbox{d}\tau \\
        = &  n|\sB(\xi^\star,1)|\cdot \int_r^\infty  e^{-\frac{\phi_t^\star + \frac{\theta \tau^2}{2}- (n+k-1)\ln(\tau^\delta)}{\delta}}\ \mbox{d}\tau,
    \end{align}\label{eq: proof-integral-ball-bound-1}\end{subequations}where the first inequality follows from a change of variables to polar coordinates and using the strong convexity of $\phi_t$ in [\ref{eq: proof-strong-convex}], and the final line by algebraic properties of logarithms.
    
    Now define
    \begin{equation}
        \varepsilon \triangleq \dfrac{\theta }{4(n+k-1)} > 0,
        \label{eq: epsilon-def}
    \end{equation}
    where the denominator is positive since $n \geq 1$ as $p$ has positive degree.
    For all $0 < \delta < \varepsilon  $,  observe
    \begin{equation}
        \tau > 1 
        \ \  \implies \ \ 
        \tau^\delta < \tau^\varepsilon  
        \ \  \mbox{and}\ \  
        \tau \leq 1
        \quad\implies\quad
        \tau^\delta \leq 1^\varepsilon,
    \end{equation}
    \ie 
    \begin{equation}
        \tau^\delta \leq \max(\tau,1)^\varepsilon,
        \quad \mbox{for all } \delta \in (0,\varepsilon).
    \end{equation}
    Whence, continuing [\ref{eq: proof-integral-ball-bound-1}], we deduce, for all $ \delta \in (0,\varepsilon),$
    \begin{subequations}
    \begin{align}
        & \int_{\bbR^n-\sB(\xi^\star,r)} e^{-\frac{\phi_t(y)}{\delta}} \|y-\xi^\star\|^k\ \mbox{d}y \\
        \leq & n|\sB(\xi^\star,1)|\cdot \int_r^\infty  e^{-\frac{\phi_t^\star + \frac{\theta \tau^2}{2} - \varepsilon(n+k-1)\ln\left(\max(\tau,1)\right)}{\delta}}\ \mbox{d}\tau.
    \end{align}\label{eq: proof-limit-laplace-argument}\end{subequations}Let $q(y)$ be the numerator inside the exponential in the integrand. Taking the limit
    \begin{equation}
        \lim_{\delta\rightarrow 0^+}  \|e^{-q}\|_{L^{\frac{1}{\delta}}([r,\infty))} 
        =  \left\| e^{-q}\right\|_{L^\infty([r,\infty))}.
    \end{equation}
    Let $\tau^\star$ be the minimizer of $q$ over $[r,\infty)$.
    If $\tau^\star > 1$, then the first order necessary condition implies, together with [\ref{eq: epsilon-def}],
    {\small 
    \begin{equation}
        0 = \theta \tau^\star -\frac{ \varepsilon(n+k-1)}{\tau^\star}
        \ \ \implies\ \
        \tau^\star = \sqrt{\dfrac{\varepsilon(n+k-1)}{\theta}}
        = \frac{1}{2},
    \end{equation}}a contradiction. Consequently, $\tau^\star \leq 1$. Since $q$ is quadratic in $\tau$ and strictly increasing on $[r,1)$, we deduce $\tau^\star = r$.
    Thus, 
    \begin{equation}
        \|e^{-q}\|_{L^\infty([r,\infty))}
        = e^{-\phi_t^\star - \frac{\theta r^2}{2}}.
        \label{eq: proof-q-delta-limit}
    \end{equation}
    Furthermore, note
    \begin{equation}
        \lim_{\delta\rightarrow0^+} \left[ n |\sB(\xi^\star,1)|\right]^\delta = 1.
        \label{eq: proof-constant-delta-limit}
    \end{equation}
    Together [\ref{eq: proof-limit-laplace-argument}], [\ref{eq: proof-q-delta-limit}], and [\ref{eq: proof-constant-delta-limit}] imply
    \begin{align}
        \lim_{\delta\rightarrow 0^+}\left[\int_{\bbR^n-\sB(\xi^\star,r)} e^{-\frac{\phi_t(y)}{\delta}} \|y-\xi^\star\|^k\ \mbox{d}y\right]^{\delta} 
        \leq  e^{-\phi_t^\star - \frac{\theta r^2}{2}}.
        \label{eq: sigma-numerator-delta-limit}
    \end{align} 
    
    \paragraph{Step 4}
     By [\ref{eq: sigma-denominator-delta-limit}] and [\ref{eq: sigma-numerator-delta-limit}] and the definition of $\sigma_\delta,$
    \begin{equation}
        \lim_{\delta\rightarrow 0^+} \left[ \int_{\bbR^n-\sB(\xi^\star,r)}\hspace*{-15pt} \sigma_\delta(y)\|y-\xi^\star\|^k\ \mbox{d}y\right]^\delta
        \leq \underbrace{ \dfrac{e^{- \phi_t^\star - \frac{\theta r^2}{2}}}{e^{ - \phi_t^\star }}}_{\triangleq \gamma}
        < 1.
    \end{equation} 
    Consequently, there is $\overline{\delta} > 0$ such that
    \begin{equation}
        \left[ \int_{\bbR^n-\sB(\xi^\star,r)}\hspace*{-15pt} \sigma_\delta(y)\|y-\xi^\star\|^k\ \mbox{d}y\right]^\delta
        \leq \dfrac{\gamma + 1}{2},
        \quad \mbox{for all\ } \delta \in (0,\overline{\delta}],
    \end{equation}
    where we note $(\gamma+1)/2 \in (\gamma,1)$, and so
    \begin{equation}
        \lim_{\delta\rightarrow 0^+}\int_{\sS}\sigma_\delta(y)\|y-\xi^\star\|^k\ \mbox{d}y
        \leq 
        \lim_{\delta \rightarrow 0^+}\left( \dfrac{\gamma+1}{2}\right)^{1/\delta}
        = 0,
    \end{equation}
    as desired.
\end{proof}

\newpage
\noindent Below we restate and prove the main theorem, which is an extension of a lemma in Section 4.5.2 of \cite{evans2010partial}.

\noindent {\bf Theorem \ref{thm: proximal-approximation}} (Proximal Approximation).
\textit{\mainresult} 

\begin{proof}  
    Let $x\in\bbR^n$ and $t > 0$ be given.
    For notational compactness, denote the HJ-prox formula by
    \begin{equation}
        \xi^\delta \triangleq \dfrac{\bbE_{y\sim  \sN(x,\delta t )}\left[{y\cdot } \exp\left(- {f}(y) / \delta \right) \right]}
    {\bbE_{y\sim  \sN(x,\delta t) }\left[ \exp\left(- {f}(y)/\delta \right) \right]},
    \quad \mbox{for all $\delta > 0$,}
    \end{equation}
    denote the proximal by
    $
        \xi^\star \triangleq \prox{tf}(x),
    $
    and note $\phi_t^\star = \phi_t(\xi^\star).$
    As argued in Lemma \ref{lemma: sigma-polynomial-vanish}, $\xi^\star$ is well-defined. We first bound $\phi_t-\phi_t^\star$ using Jensen's inequality (Step 1).
    Second, we show $\phi_t(\xi^\delta) \rightarrow \phi_t(\xi^\star)$ (Step 2). The strong convexity of  $\phi_t$  enables us to establish the desired limit (Step 3). 
  
    \paragraph{Step 1} 
    Note $\xi^\delta$ can be rewritten via
    \begin{equation}
        \xi^\delta = \left[\int_{\bbR^n} e^{-\frac{\phi(y)}{\delta}}\mbox{d}y \right]^{-1}   \int_{\bbR^n} y \cdot e^{ -\frac{\phi_t(y)}{\delta}}\ \mbox{d}y.
    \end{equation}
    Using $\sigma_\delta$,  the estimate can   be more concisely written via
    \begin{equation}
        \xi^\delta = \int_{\bbR^n} \sigma_\delta(y) y \ \mbox{d}y
        = \bbE_{y\sim\bbP_{\sigma_\delta}}\left[y\right],
    \end{equation}
    where the expectation holds by utilizing the fact [\ref{eq: proof-sigma-pdf}] shows $\sigma_\delta$ defines a probability density. Thus,   Jensen's inequality may be applied to reveal  
    \begin{equation}
        0 \leq \phi_t^\star \leq \phi_t(\xi^\delta)
        = \phi_t\left( \bbE_{y\sim\sigma_\delta} [y]\right)
        \leq \bbE_{y\sim\sigma_\delta}\left[ \phi_t(y)\right].
    \end{equation}
    In integral form, we may subtract $\phi_t^\star$ to write
    \begin{equation}
        0 \leq \phi_t(\xi^\delta) - \phi_t^\star
        \leq \int_{\bbR^n}\sigma_\delta(y)[ \phi_t(y) - \phi_t^\star]\ \mbox{d}y.
        \label{eq: proof-Jensen}
    \end{equation}

    \paragraph{Step 2}
    Let $\varepsilon > 0$ be given. To deduce $\phi_t(\xi^\delta)\rightarrow \phi_t^\star$, we verify there is $\delta^\star > 0$ such that
    \begin{equation}
        |\phi_t(\xi^\delta) - \phi_t^\star| \leq \varepsilon,
        \quad \mbox{for all}\ \delta \in (0,\delta^\star].
        \label{eq: proof-function-convergence}
    \end{equation}
    By [\ref{eq: proof-Jensen}], the relation [\ref{eq: proof-function-convergence}] holds
    if there is such a $\delta^\star$ that
    \begin{equation}
        \int_{\bbR^n}\sigma_\delta(y)[ \phi_t(y) - \phi_t^\star]\ \mbox{d}y
        \leq \varepsilon,
        \quad \mbox{for all }\ \delta \in (0,\delta^\star].
        \label{eq: proof-obj-conv}
    \end{equation}
    We verify this   by splitting the integral into two parts. Since $f$ is either $L$-Lipschitz or $L$-smooth, there is a quadratic polynomial $p\colon\bbR\rightarrow\bbR$ with nonnegative coefficients such that
    \begin{equation}
        \phi_t(y) \leq p(\|y-\xi^\star\|),
        \quad \mbox{for all $y\in\bbR^n$,}
    \end{equation}
    and $p(0) = \phi_t^\star$. Thus, by the intermediate value theorem, we may fix $r \in (0,1)$ sufficiently small to ensure  
    \begin{equation}
       p(r) -\phi_t^\star \leq \dfrac{\varepsilon}{2}.
    \end{equation} 
    This implies
    \begin{subequations}
        \begin{align}
        \label{eq: proof-phi-bound}
        \phi_t(y) - \phi_t^\star 
        & \leq p(\|y-\xi^\star\|) -\phi_t^\star \\
        &\leq \dfrac{\varepsilon}{2} ,
        \quad \mbox{for all} \ y\in\sB(\xi^\star,r).
    \end{align}\end{subequations}%
    Thus, integrating over the ball $\sB(\xi,r)$ reveals
    \begin{subequations}
        \begin{align}
            A & \triangleq \int_{\sB(\xi^\star,r)} \sigma_\delta(y) [\phi_t(y)-\phi_t^\star] \ \mbox{d}y \\
            & \leq \int_{\sB(\xi,r)} \sigma_\delta(y)\cdot \dfrac{\varepsilon}{2} \ \mbox{d}y \\
            & \leq \dfrac{\varepsilon}{2}\cdot  \int_{\bbR^n} \sigma_\delta(y)  \ \mbox{d}y \\  
            &  = \dfrac{\varepsilon}{2},
        \end{align}\label{eq: proof-A}\end{subequations}where the second inequality follows from [\ref{eq: proof-sigma-pdf}]. Next we  integrate over the rest of $\bbR^n$. Define
    \begin{subequations}
        \begin{align}
            B_\delta 
            & \triangleq \int_{\bbR^n-\sB(\xi^\star,r)} \sigma_\delta(y) [\phi_t(y)-\phi_t^\star]\ \mbox{d}y \\ 
            & \leq \int_{\bbR^n-\sB(\xi^\star,r)} \sigma_\delta(y)\cdot p(\|y-\xi^\star\|)\ \mbox{d}y.
        \end{align}\end{subequations}We may apply Lemma \ref{lemma: sigma-polynomial-vanish} to deduce there is $\omega > 0$ such that 
    \begin{equation}
        B_\delta \leq \dfrac{\varepsilon}{2},
        \quad \mbox{for all } \delta \in (0,\omega].
        \label{eq: proof-B}
    \end{equation} 
    Consequently, [\ref{eq: proof-A}] and [\ref{eq: proof-B}] together imply
    \begin{subequations}
        \begin{align}
            \int_{\bbR^n} \sigma_\delta(y) [ \phi_t(y) - \phi_t^\star] \ \mbox{d}y 
            & = A + B_\delta \\
            & \leq \dfrac{\varepsilon}{2} +\dfrac{\varepsilon}{2}  \\
            & \leq \varepsilon,
            \quad \mbox{for all $\delta \in (0, \omega].$}
        \end{align}
    \end{subequations}
    Hence [\ref{eq: proof-obj-conv}] holds, taking $\delta^\star = \omega,$ \ie $\phi_t(\xi^\delta) \rightarrow \phi_t^\star$ as $\delta\rightarrow 0^+$.

    \paragraph{Step 3}Let $\overline{\varepsilon} > 0$. It suffices to show there is $\overline{\delta} > 0$ such that
    \begin{equation}
        \| \xi^\delta - \xi^\star\| \leq \overline{\varepsilon},
        \quad \mbox{for all $\delta \in (0, \overline{\delta}].$}
        \label{eq: proof-01}
    \end{equation} 
    Define
    \begin{equation}
        \sS \triangleq \left\lbrace  z : \|z-\xi^\star\| \geq  \overline{\varepsilon}
        \right\rbrace 
    \end{equation}
    and note, by the strong convexity  of $\phi_t$ (\eg see [\ref{eq: proof-strong-convex}]),
    \begin{equation}
        \phi_t(z) \geq \phi_t^\star + \dfrac{\theta\overline{\varepsilon}^2}{2},
        \quad \mbox{for all $z\in\sS$.}
    \end{equation} 
    By Step 2, there is $\mu  > 0$ such that
    \begin{equation}
        \phi_t(\xi^\delta)
        \leq \phi_t^\star  + \dfrac{\theta \overline{\overline{\varepsilon}}^2}{4}, 
        \quad \mbox{for all $\delta \in (0, \mu ].$}
    \end{equation}
    Thus, $\xi^\delta \notin \sS$, for all $\delta \in (0, \mu]$,
    \ie  (\ref{eq: proof-01}) holds, taking $\overline{\delta} = \mu.$ This completes the proof. 
\end{proof}  
\end{document}